\documentclass[12pt]{amsart}
\usepackage{amsmath, amssymb, latexsym, amsthm}
\newif\iflogic
\logicfalse   %%%Use mathematical notation and include appendix

\newcommand{\arx}[1]{\texttt{http://arxiv.org/abs/#1}}
\newcommand{\depth}{\op{Depth^+}}

\newcommand{\cB}{\mathcal{B}}

\newcommand{\cU}{\mathcal{U}}
\newcommand{\fI}{\mathfrak{I}}
\newcommand{\fJ}{\mathfrak{J}}
\newcommand{\fF}{\mathfrak{F}}
\newcommand{\fD}{\mathfrak{D}}
\newcommand{\fB}{\mathfrak{B}}
\newcommand{\Dfin}{\mathfrak{D}_\mathrm{fin}}
\newcommand{\fX}{\mathfrak{X}}
\renewcommand{\ss}{\hspace{0.07cm}}
\renewcommand{\mid}[5]{[{#4}\ss{#1}\ss{#3}\ss{#2}\ss{#5}]}
\newcommand{\domi}[3]{[{#1}\ss{#2}\ss{#3}]}

\newcommand{\emp}{excluded middle property}

\newcommand{\wemp}{weak \emp{}}

\newcommand{\op}{\operatorname}
\newcommand{\maxfin}{\op{maxfin}}

\newcommand{\Cal}{\mathcal}

\newcommand{\bbC}{\mathbb{C}}

\newcommand{\cF}{{\Cal F}}

\newcommand{\M}{{\Cal M}}
\newcommand{\Null}{{\Cal N}}
\newcommand{\N}{\iflogic\omega\else\mathbb{N}\fi}
\newcommand{\NN}{{{}^{\N}\N}}

\renewcommand{\inf}{\iflogic{[\omega]^\omega}\else{P_\oo(\N)}\fi}

\newcommand{\Union}{\bigcup}

\newcommand{\Impl}{\Rightarrow}
\long\def\forget#1\forgotten{}
\renewcommand{\b}{{\mathfrak b}}
\renewcommand{\u}{{\mathfrak u}}
\renewcommand{\r}{{\mathfrak r}}

\newcommand{\kwt}{\kappa_{\w\tau}}
\newcommand{\fx}{\mathfrak{x}}

\renewcommand{\d}{{\mathfrak d}}

\newcommand{\g}{\mathfrak{g}}
\renewcommand{\i}{\item}

\newcommand{\oo}{\infty}
\newcommand{\p}{{\mathfrak p}}
\newcommand{\s}{\mathfrak{s}}
\renewcommand{\a}{\mathfrak{a}}
\newcommand{\e}{\mathfrak{e}}
\newcommand{\w}{\omega}

\newcommand{\nin}{\not\in}

\newcommand{\sbst}{\subseteq}
\newcommand{\spst}{\supseteq}
\newcommand{\sm}{\setminus}

\newcommand{\as}{\subseteq^*}%{\let\proclaim\relax}
\renewcommand{\pi}{pseudo-intersection}
\renewcommand{\|}{\restriction}

\newcommand{\<}{\langle}
\renewcommand{\>}{\rangle}

\newcommand{\cov}{{\sf cov}}
\newcommand{\add}{{\sf add}}
\newcommand{\cof}{{\sf cof}}
\newcommand{\cf}{{\sf cf}}
\newcommand{\non}{{\sf non}}

\newcommand{\impl}{\to}
\renewcommand{\t}{\mathfrak{t}}

\newtheorem{thm}{Theorem}%[section]

\newtheorem{prob}[thm]{Problem}
\newtheorem{lem}[thm]{Lemma}
\newtheorem{cor}[thm]{Corollary}
\theoremstyle{definition}
\newtheorem{defn}[thm]{Definition}
\newtheorem{rem}[thm]{Remark}
\newcommand{\be}{\begin{enumerate}}
\newcommand{\ee}{\end{enumerate}}
\newcommand{\bi}{\begin{itemize}}
\newcommand{\ei}{\end{itemize}}

\title[Combinatorial notions of smallness and additivity]
{Critical cardinalities and additivity properties
of combinatorial notions of smallness}

\begin{document}
\author{Saharon Shelah}
\thanks{The research of the first author is partially
supported by The Israel Science Foundation founded
by the Israel Academy of Sciences and Humanities. Publication 768.}
\address{
Institute of Mathematics, Hebrew University of Jerusalem,
Givat Ram, 91904 Jerusalem, Israel,
and Mathematics Department, Rutgers University,
New Brunswick, NJ \ 08903, U.S.A.
}
\email{shelah@math.huji.ac.il}
%\and
\author{Boaz Tsaban}
\thanks{This paper constitutes a part of the second author's doctoral dissertation at
Bar-Ilan University.}
\address{Department of Mathematics and Computer Science, Bar-Ilan University,
Ramat-Gan 52900, Israel}
\email{tsaban@macs.biu.ac.il}

\begin{abstract}
Motivated by the minimal tower problem,
an earlier work studied diagonalizations of covers
where the covers are related to linear
quasiorders ($\tau$-covers).
We deal with two types of combinatorial questions which arise
from this study.
\be
\i Two new cardinals introduced in
the topological study are expressed
in terms of well known cardinals characteristics of the continuum.
\i We study the additivity numbers of the combinatorial
notions corresponding to the topological diagonalization
notions.
\ee
This gives new insights on the structure of the eventual dominance
ordering on the Baire space, the almost inclusion ordering on the
Rothberger space, and the interactions between them.
\end{abstract}


\maketitle
%\markboth{\sc Boaz Tsaban}{\sc Selection principles and the minimal tower problem}

\section{Introduction and overview}

\iflogic
Let $\N$ denote the set of natural numbers.
\else
\fi
We work with two spaces which carry an interesting combinatorial structure:
The \emph{Baire space} $\NN$ with  eventual dominance $\le^*$
($f\le^* g$ if $f(n)\le g(n)$ for all but finitely many $n$),
and the \emph{Rothberger space} $\inf=\{A\sbst\N : A\mbox{ is infinite}\}$
with $\as$ ($A\as B$ if $A\sm B$ is finite).
We write $A\subset^* B$ if $A\as B$ and $B\not\as A$.

A subset $X$ of $\NN$ is \emph{unbounded} if it is unbounded with
respect to $\le^*$.
$X$ is \emph{dominating} if it is cofinal in $\NN$ with respect to $\le^*$.
$\b$ is the minimal size of an unbounded subset of $\NN$, and $\d$ is the minimal size
of a dominating subset of $\NN$.

An infinite set $A\sbst\N$ is a \emph{pseudo-intersection}
of a family $\cF\sbst\inf$
if for each $B\in \cF$, $A\as B$.
A family $\cF\sbst\inf$ is a \emph{tower} if it is linearly quasiordered
by $\as$, and it has no \pi{}.
$\t$ is the minimal size of a tower.
A family $\cF\sbst\inf$ is \emph{centered} if the intersection of
each (nonempty) finite subfamily of $\cF$ is infinite.
$\p$ is the minimal size of a centered family which has no \pi{}.
A family $\cF\sbst\inf$ is \emph{splitting} if for each infinite $A\sbst\N$
there exists $S\in\cF$ which \emph{splits} $A$, that is,
such that the sets $A\cap S$ and $A\sm S$ are infinite.
$\s$ is the minimal size of a splitting family.

Let $\mathfrak{c}=2^{\aleph_0}$.
The following relations,
where an arrow means $\le$,  are well-known \cite{HBK}:
$$\begin{matrix}
 & & & & & & \b\\
 & & & & & \nearrow & & \searrow\\
 \aleph_1 & \impl & \p & \impl & \t & & & & \d & \impl \mathfrak{c}\\
 & & & & & \searrow & & \nearrow\\
 & & & & & & \s
\end{matrix}$$
No pair of cardinals in this diagram is provably equal, except
perhaps $\p$ and $\t$. The \emph{Minimal Tower problem},
which asks whether it is provable that $\p=\t$,
is one of the most important
problems in infinite combinatorics, and it
goes back to Rothberger (see, e.g., \cite{ROTH2}).

\subsection*{New cardinals}
In \cite{tsaban1}, topological notions related
to $\p$ and $\t$ were compared. In \cite{tautau} the
topological notion related to $\t$ (called \emph{$\tau$-covers})
was studied in a wider
context. This study led back to several new combinatorial
questions, one of which related to the
minimal tower problem.

\begin{defn}
For a family $\cF\sbst\inf$ and an infinite $A\sbst\N$, define
$\cF\| A = \{B\cap A : B\in \cF\}$.
If all sets in $\cF\| A$ are infinite, we
say that $\cF\| A$ is a
\emph{large restriction} of $\cF$.
Let $\kwt$ be the minimal cardinality of a \emph{centered} family $\cF\sbst\inf$
such that there exists no infinite $A\sbst\N$ such that the restriction
$\cF\| A$ is large and linearly quasiordered by $\as$.
\end{defn}

It is not difficult to see that $\p=\min\{\kwt,\t\}$ \cite{tautau}.
In Section \ref{kwt} we show that in fact, $\p=\kwt$.
This existence of a centered family with no large
linearly quasiordered restriction shows that
$\p$ is combinatorially ``larger'' than asserted in its
original definition,
and suggests an additional evidence to the difficulty
of separating $\p$ from the combinatorially ``larger''
cardinal $\t$: Now the consistency of $\kwt<\t$ must be established
in order to solve the Minimal Tower problem in the negative.

\begin{defn}
For functions $f,g\in\NN$, and a binary relation $R$
on $\N$, define a subset $\domi{f}{R}{g}$
of $\N$ by:
$$\domi{f}{R}{g} = \{ n : f(n)Rg(n)\}.$$
Next, For functions $f,g,h\in\NN$, and binary relations $R,S$
on $\N$, define $\mid{R}{S}{g}{h}{f}\sbst\N$ by:
$$\mid{R}{S}{g}{f}{h}=\domi{f}{R}{g}\cap\domi{g}{S}{h}
=\{n : f(n)R g(n)\mbox{ and }g(n)S h(n)\}.$$
For a subset $X$ of $\NN$ and $g\in\NN$, we say that
$g$ \emph{avoids middles} in $X$ with respect to $\<R,S\>$ if:
\be
\i for each $f\in X$, the set $\domi{f}{R}{g}$ is infinite;
\i for all $f,h\in X$ at least one of the sets
$\mid{R}{S}{g}{f}{h}$ and $\mid{R}{S}{g}{h}{f}$ is finite.
\ee
$X$ satisfies the \emph{excluded middle} property with respect to $\<R,S\>$
if there exists $g\in\NN$ which avoids middles in $X$
with respect to $\<R,S\>$.
$\fx_{R,S}$ is the minimal size of a subset $X$ of $\NN$
which does not satisfy the \emp{} with respect to
$\<R,S\>$.
\end{defn}
The cardinal $\fx=\fx_{<,\le}$ was defined in \cite{tautau}.
In section \ref{excluded} we express
all of the four cardinals
$\fx_{\le,\le},\fx_{<,\le},\fx_{\le,<}$, and $\fx_{<,<}$ in terms of
well-known cardinals. This solves several problems raised in \cite{tautau}.

\subsection*{Additivity of combinatorial notions of smallness}
For a finite subset $F$ of $\NN$, define $\max(F)\in\NN$ by
$\max(F)(n)=\max\{f(n) : f\in F\}$ for each $n$.
A subset $Y$ of $\NN$ is \emph{finitely-dominating}
if the collection
$$\maxfin(Y):=\{\max(F) : F\mbox{ is a finite subset of }Y\}$$
is dominating.

We will use the following notations:
\bi
\i[$\fB$]: The collection of all bounded subsets of $\NN$,
\i[$\fX$]: The collection of all subsets of $\NN$ which satisfy the
excluded middle property with respect to $\<<,\le\>$,
\i[$\Dfin$]: The collection of all subsets of $\NN$ which are not
finitely dominating; and
\i[$\fD$]: The collection of all subsets of $\NN$ which are not dominating.
\ei
Thus $\fB\sbst\fX\sbst\Dfin\sbst\fD$.
The classes $\fB$, $\fX$, $\Dfin$, and $\fD$ are used to characterize
certain topological diagonalization properties \cite{CBC, huremen1, tautau}.

Following \cite{JUBAR}, we define the \emph{additivity number}
for classes $\fI\sbst\fJ\sbst P(\NN)$ with $\cup\fI\nin\fJ$ by
$$\add(\fI,\fJ)=\min\{|\fF| : \fF\sbst\fI\mbox{ and }\cup\fF\nin\fJ\},$$
and write $\add(\fJ)=\add(\fJ,\fJ)$.
If $\fI$ contains all singletons, then $\add(\fI,\fJ)\le\non(\fJ)$,
where $\non(\fJ)=\min\{|J| : J\sbst\NN\mbox{ and }J\nin\fJ\}$
(thus $\non(\fB)=\b$, $\non(\fD)=\non(\Dfin)=\d$,
and $\non(\fX)=\fx$.)

For $\fI,\fJ\in\{\fB,\fX,\Dfin,\fD\}$, the cardinals
$\add(\fI,\fJ)$ bound from below the additivity
numbers of the corresponding topological diagonalizations.
In section \ref{addsection} we express $\add(\fI,\fJ)$
for almost all $\fI,\fJ\in\{\fB,\fX,\Dfin,\fD\}$ in terms
of well known cardinal characteristics of the continuum.
In two cases for which this is not done, we give
consistency results.

\section{The cardinal $\kwt$}\label{kwt}
For our purposes, a \emph{filter} on a boolean subalgebra $\cB$
of $P(\N)$ is a family $\cU\sbst \cB$ which is closed under
taking supersets in $\cB$ and finite intersections, and does not
contain finite sets as elements.

\begin{thm}
$\p=\kwt$.
\end{thm}
\begin{proof}
Let $\cF\sbst\inf$ be a centered family of size $\p$ which has no \pi{}.
Let $\cB$ be the boolean subalgebra of $P(\N)$ generated by $\cF$.
Then $|\cB|=\p$.
Let $\cU\sbst \cB$ be a filter of $\cB$ containing $\cF$.
As $\cU$ does not contain finite sets as elements,
$\cU$ is centered. Moreover, $|\cU|=\p$, and it has
no \pi{}.

Towards a contradiction, assume that $\p<\kwt$.
Then there exists an infinite $A\sbst\N$ such that the restriction
$\cU\| A$ is large, and is linearly quasiordered by $\as$.
Fix any element $D_0\cap A\in \cU\|A$.
As $\cU\| A$ does not have a \pi{},
there exist:
\be
\i An element $D_1\cap A\in \cU\| A$
such that $D_1\cap A\subset^* D_0\cap A$;
and
\i An element $D_2\cap A\in \cU\| A$ such that $D_2\cap A\subset^* D_1\cap A$.
\ee
Then the sets $(D_2\cup(D_0\sm D_1))\cap A$ and
$D_1\cap A$ (which are elements of $\cU\| A$)
contain the infinite sets
$(D_0\cap A)\sm (D_1\cap A)$ and $(D_1\cap A)\sm (D_2\cap A)$,
respectively, and thus are not $\as$-comparable, a contradiction.
\end{proof}

A closely related problem from \cite{tautau} remains open.

\begin{defn}
A family $Y\subseteq\inf$ is \emph{linearly refinable}
if for each $y\in Y$ there exists an infinite subset
$\hat y\subseteq y$ such that the family $\hat Y = \{\hat y : y\in Y\}$ is
linearly $\as$-quasiordered.
$\p^*$ is the minimal size of a centered family
in $\inf$ which is not linearly refineable.
\end{defn}

Again it is easy to see that $\p = \min\{\p^*,\t\}$.
Thus, a solution of the following problem may shed more light on the Minimal Tower problem.

\begin{prob}
Does $\p=\p^*$?
\end{prob}

\section{The excluded middle property}\label{excluded}

\begin{lem}\label{easylemma}
$\b\le\fx_{\le,\le}\le\fx_{\le,<}\le\fx_{<,\le}\le\fx_{<,<}\le\d$.
\end{lem}
\begin{proof}
The inequalities $\fx_{\le,\le}\le\fx_{\le,<}$ and $\fx_{<,\le}\le\fx_{<,<}$
are immediate from the definitions. We will prove the other inequalities.

Assume that $Y$ is a bounded subset of $\NN$. Let $g\in\NN$
bound $Y$. Then $g$ avoids middles in $Y$
with respect to $\<\le,\le\>$.
This shows that $\b\le\fx_{\le,\le}$.

Next, consider a subset $Y$ of $\NN$ which satisfies the
\emp{} with respect to $\<<,<\>$, and let $g$ witness that.
Then $g$ witnesses that $Y$ is not dominating.
Thus $\fx_{<,<}\le\d$.

It remains to show that $\fx_{\le,<}\le\fx_{<,\le}$.
Assume that $Y\sbst \NN$ satisfies the \emp{} with respect to $\<\le,<\>$,
and let $g\in\NN$ avoid middles in $Y$ with respect to
$\<\le,<\>$.
Define $\tilde g\in\NN$
such that $\tilde g(n)=g(n)+1$ for each $n$.
For each $f,h\in Y$ we have that $\domi{f}{\le}{g} = \domi{f}{<}{\tilde g}$,
and $\mid{\le}{<}{g}{f}{h}=\mid{<}{\le}{\tilde g}{f}{h}$.
Therefore,
$\tilde g$ avoids middles in $Y$ with respect to $\<<,\le\>$.
\end{proof}

\begin{thm}\label{equalsb}
$\fx_{\le,\le}=\fx_{\le,<}=\b$.
\end{thm}
\begin{proof}
By Lemma \ref{easylemma},
it is enough to show that $\fx_{\le,<}\le\b$.
Let $\<b_\alpha : \alpha<\b\>$ be an unbounded subset of $\NN$.
For each $\alpha<\b$ define $b^0_\alpha,b^1_\alpha\in\NN$ by
$$\begin{cases}
b^0_\alpha(2n) & = b_\alpha(n)\\
b^0_\alpha(2n+1) & = 0
\end{cases};
\quad
\begin{cases}
b^1_\alpha(2n) & = 0\\
b^1_\alpha(2n+1) & = b_\alpha(n)
\end{cases}$$
%\begin{eqnarray*}
%
%
%\end{eqnarray*}
for each $n\in\N$,
and set $Y = \{b^0_\alpha, b^1_\alpha : \alpha<\b\}$. Then $|Y|=\b$.
We will show that $Y$ does not satisfy the
excluded middle property with respect to $\<\le,<\>$.
For each $g\in\NN$, let $\alpha<\b$ be such that
$\max\{g(2n),g(2n+1)\}<b_\alpha(n)$
for infinitely many $n$.
Then:
\begin{eqnarray*}
\mid{\le}{<}{g}{b^0_\alpha}{b^1_\alpha}
& =     & \{ n : b^0_\alpha(n)\le g(n)<b^1_\alpha(n)\}\\
& \spst & \{2n+1 : 0\le g(2n+1)<b_\alpha(n)\}
\end{eqnarray*}
is an infinite set.
Similarly,
$\mid{\le}{<}{g}{b^1_\alpha}{b^0_\alpha}\spst \{2n : 0\le g(2n)<b_\alpha(n)\}$
is also infinite. That is, $g$ does not avoid middles in $Y$
with respect to $\<\le,<\>$.
\end{proof}

\begin{lem}\label{s-le}
$\s\le\fx_{<,\le}$.
\end{lem}
\begin{proof}
Assume that $Y\sbst\NN$ is such that $|Y|<\s$.
%%Trying a shorter proof
Let $\cF\sbst P(\N)$ be the family of all sets
of the form $\domi{f}{<}{h}$, %$\{n : f(n)<h(n)\}$,
where $f,h\in Y$.
$|\cF|<\s$, thus there exists an infinite subset $A$ of $\N$
such that for each $X\in \cF$, either $A\cap X$ is finite, or $A\sm X$
is finite.
As $|Y|<\s\le\d$, there exists $g\in\NN$ such that for each
$f\in Y$, $g\| A\not\le^* f\| A$.
(In particular, $\domi{f}{<}{g}$ is infinite for each $f\in Y$.)
We may assume that for $n\nin A$, $g(n)=0$.

Consider any set $\domi{f}{<}{h}\in \cF$.
If $A\cap\domi{f}{<}{h}$ is finite, then the set
\begin{eqnarray*}
\mid{<}{\le}{g}{f}{h}
%& =     & \{ n : f(n)<g(n)\le h(n) \}\\
& \sbst & \{ n : 0<g(n), f(n) < h(n) \}\\
& \sbst & \{ n\in A : f(n) < h(n) \} = A\cap\domi{f}{<}{h}
\end{eqnarray*}
is finite.
Otherwise, $A\sm\domi{f}{<}{h}$ is finite, so we get similarly that
\begin{eqnarray*}
\mid{<}{\le}{g}{h}{f}
& \sbst & \{ n\in A : h(n) < f(n) \}\\
& \sbst & \{ n\in A : h(n) \le f(n) \} = A\sm\domi{f}{<}{h}
\end{eqnarray*}
is finite. Thus $Y$ satisfies the \emp{} with respect to
$\<<,\le\>$.
\end{proof}

\begin{thm}\label{sbthm}
$\fx_{<,\le}=\fx_{<,<}=\max\{\s,\b\}$.
\end{thm}
\begin{proof}
By Lemmas \ref{easylemma} and \ref{s-le}, we have that
$\max\{\s,\b\}\le\fx_{<,\le}\le\fx_{<,<}$.
We will prove that $\fx_{<,<}\le\max\{\s,\b\}$.
The argument is an extension of the proof of
Theorem \ref{equalsb}.

Let
$\b^*$
be the minimal size of a subset $B$ of $\NN$ such that
$B$ is unbounded on each infinite subset of $\N$.
According to \cite{HBK}, $\b=\b^*$.
Thus there exists a subset
$B=\<b_\alpha : \alpha<\b\>$ of $\NN$
such that $B$ is increasing with respect to
$\le^*$ and unbounded on each infinite subset of $\N$.
Let $\mathcal{S} = \<S_\alpha : \alpha<\s\>\sbst\inf$ be a splitting family.
For each $\alpha<\s$ and $\beta<\b$ define
$b^0_{\alpha,\beta},b^1_{\alpha,\beta}\in\NN$ by:
$$b^0_{\alpha,\beta}(n) =
\begin{cases}
b_\beta(n) & n\in S_\alpha\\
0          & n\nin S_\alpha
\end{cases};
\quad
b^1_{\alpha,\beta}(n) =
\begin{cases}
0          & n\in S_\alpha\\
b_\beta(n) & n\nin S_\alpha
\end{cases}$$
and set $Y = \{b^i_{\alpha,\beta} : i<2,\alpha<\s,\beta<\b\}$.
Then $|Y|=2\cdot\s\cdot\b=\max\{\s,\b\}$.
We will show that $Y$ does not satisfy the \emp{} with respect to $\<<,<\>$.
Assume that $g\in\NN$ avoids middles in $Y$ with respect to
$\<<,<\>$. Then the set $A=\domi{0}{<}{g}$ is infinite; thus there
exists $\alpha<\s$ such that the sets $A\cap S_\alpha$ and $A\sm S_\alpha$
are infinite.
Pick $\gamma<\b$ such
that $b_\gamma\|A\cap S_\alpha\not\le^* g\|A\cap S_\alpha$,
and $\beta>\gamma$ such that $b_\beta\|A\sm S_\alpha\not\le^* g\|A\sm S_\alpha$.
%, and set $\gamma=\max\{\alpha,\beta\}$.
Then:
\begin{eqnarray*}
\mid{<}{<}{g}{b^0_{\alpha,\beta}}{b^1_{\alpha,\beta}}
& \spst &
\{ n\in A\sm S_\alpha : b^0_{\alpha,\beta}(n)<g(n)<b^1_{\alpha,\beta}(n)\}\\
& = & \{ n\in A\sm S_\alpha : 0 <g(n) <b_\beta(n)\}\\
& = & \{ n\in A\sm S_\alpha : g(n) <b_\beta(n)\}
\end{eqnarray*}
is an infinite set.
Similarly, the set
\begin{eqnarray*}
\mid{<}{<}{g}{b^1_{\alpha,\beta}}{b^0_{\alpha,\beta}}
& \spst &
\{ n\in A\cap S_\alpha : b^1_{\alpha,\beta}(n)<g(n)<b^0_{\alpha,\beta}(n)\}\\
& = & \{ n\in A\cap S_\alpha : 0 <g(n) <b_\beta(n)\}\\
& = & \{ n\in A\cap S_\alpha : g(n) <b_\beta(n)\}
\end{eqnarray*}
is also infinite, because $b_\gamma\le^* b_\beta$;
a contradiction.
\end{proof}

\begin{rem}
The cardinal $\max\{\s,\b\}$ is also equal to the
\emph{finitely splitting number} $\mathfrak{fs}$ studied in
\cite{KaWe}.
\end{rem}

Several variations of the \emp{} are studied in the appendix to the online
version of this paper \cite{online}.

\section{Additivity of combinatorial properties}\label{addsection}

The additivity number $\add(\fI,\fJ)$ is monotone
decreasing in the first coordinate and increasing in the second.
Our task in this section is to determine, when possible, the
cardinals in the following diagram in terms of the usual cardinal
characteristics $\b$, $\d$, etc. (In this diagram, an
arrow means $\le$.)

$$\begin{matrix}
\add(\fD,\fD)&\impl&\add(\Dfin,\fD)   &\impl&\add(\fX,\fD)&\impl&\add(\fB,\fD)\\
           &     &\uparrow         &     & \uparrow  &     & \uparrow\\
           &     &\add(\Dfin,\Dfin)&\impl&\add(\fX,\Dfin)&\impl&\add(\fB,\Dfin)\\
           &     &                 &     & \uparrow  &     & \uparrow\\
           &     &                 &     &\add(\fX,\fX)&\impl&\add(\fB,\fX)\\
           &     &                 &     &           &     & \uparrow\\
           &     &                 &     &           &     &\add(\fB,\fB)
\end{matrix}$$

\subsection{Results in ZFC}
\begin{thm}
The following equalities hold:
\be
\i $\add(\fB,\Dfin)=\add(\fB,\fD)=\d$,
\i $\add(\Dfin,\Dfin)=\add(\fX,\fX)=\add(\fX,\Dfin)=2$; and
\i $\add(\fD,\fD)=\add(\fB,\fB)=\add(\fB,\fX)=\b$.
\ee
\end{thm}
\begin{proof}
(1) As $\non(\fD)=\d$, it is enough to show that $\add(\fB,\Dfin)\ge\d$.
%Observe that for each $Y\in\fB$, $\maxfin(Y)\in\fB$ as well.
Assume that $|I|<\d$, and that $Y=\Union_{i\in I} Y_i$ where
each $Y_i$ is bounded. For each $i\in I$ let $g_i$ bound
$Y_i$. As $|I|<\d$, the family $\maxfin(\{g_i : i\in I\})$ is
not dominating; let $h$ be a witness for that.
For each finite $F\sbst Y$, let $\tilde I$ be a finite subset
of $I$ such that $F\sbst\Union_{i\in\tilde I} Y_i$.
Then $\max(F)\le^*\max(\{g_i : i\in\tilde I\})\not\ge^* h$.
Thus $\max(F)\not\ge^* h$, so $Y\in\Dfin$.

(2) It is enough to show that $\add(\fX,\Dfin)=2$.
Thus, let
\begin{eqnarray*}
Y_0 & = & \{f\in\NN : (\forall n)f(2n)=0\mbox{ and }f(2n+1)\ge 1\}\\
Y_1 & = & \{f\in\NN : (\forall n)f(2n)\ge 1\mbox{ and }f(2n+1)=0\}
\end{eqnarray*}
Then the constant function $g\equiv 1$ witnesses that $Y_0,Y_1\in\fX$,
but $Y_0\cup Y_1$ is $2$-dominating, and in particular finitely dominating.

(3) It is folklore that $\add(\fD,\fD)=\add(\fB,\fB)=\b$ -- see, e.g.,
\cite{AddQuad} for a proof.
It remains to show that $\add(\fB,\fX)\le\b$.
Let $B$ be a subset of $\NN$ which is unbounded on each
infinite subset of $\N$, and such that $|B|=\b$.
For each $f\in B$ let $Y_f = \{g\in\NN : g\le^* f\}$. (Thus
each $Y_f$ is bounded.)
We claim that $Y=\Union_{f\in B}Y_f\nin\fX$.
To this end, consider any function $g\in\NN$ which claims to
witness that $Y\in\fX$. In particular, $\domi{0}{<}{g}$ must be
infinite. Choose $f\in B$ such that
$f\|\domi{0}{<}{g}\not\le^* g\|\domi{0}{<}{g}$, that
is, $\mid{<}{<}{g}{0}{f}$ is infinite.
Let $A_0,A_1$ be a partition of $\mid{<}{<}{g}{0}{f}$ into two infinite
sets, and define $f_0\in Y_f$ by $f_0(n) = g(n)$ when $n\in A_0$
and $0$ otherwise; similarly define
$f_1\in Y_f$ by $f_1(n) = g(n)$ when $n\in A_1$
and $0$ otherwise. Then $f_0,f_1\in Y$, but both of the sets
$\mid{<}{\le}{g}{f_0}{f_1}$ and
$\mid{<}{\le}{g}{f_1}{f_0}$ are infinite.
\end{proof}

\subsection{Consistency results}
The only cases which we have not solved yet are
$\add(\Dfin,\fD)$ and $\add(\fX,\fD)$.
In \cite{AddQuad} it was proved that $\b\le\add(\Dfin,\fD)$.
In Theorem 2.2 of \cite{Mildenberger} it is (implicitly) proved that $\g\le\add(\Dfin,\fD)$.
Thus
$$\max\{\b,\g\}\le\add(\Dfin,\fD)\le\add(\fX,\fD)\le\d.$$
Moreover, for any $\fI\sbst\fJ$, $\cf(\add(\fI,\fJ))\ge\add(\fJ)$,
and therefore
$$\cf(\add(\Dfin,\fD)),\cf(\add(\fX,\fD))\ge\add(\fD,\fD)=\b.$$

The notion of ultrafilter will be used to obtain upper bounds
on $\add(\Dfin,\fD)$ and $\add(\fX,\fD)$. A family
$\cU\sbst\inf$ is a \emph{nonprincipal ultrafilter} if
it is closed under taking supersets and finite intersections,
and cannot be extended, that is, for each infinite $A\sbst\N$,
either $A\in\cU$ or $\N\sm A\in\cU$.
Consequently, a linear quasiorder $\le_\cU$ can be defined
on $\NN$ by
$$f\le_\cU g\mbox{\qquad if\qquad}\domi{f}{\le}{g}\in\cU.$$
The \emph{cofinality} of the reduced product $\NN/\cU$ is
the minimal size of a subset $C$ of $\NN$ which is
cofinal in $\NN$ with respect to $\le_\cU$.

\begin{thm}\label{filteraddDfinD}
For each cardinal number $\kappa$, the following are equivalent:
\be
\i $\kappa<\add(\Dfin,\fD)$;
\i For each $\kappa$-sequence $\<(g_\alpha,\cU_\alpha) : \alpha<\kappa\>$
with each $\cU_\alpha$ an %nonprincipal
ultrafilter on $\N$ and each $g_\alpha\in\NN$
there exists $g\in\NN$ such that for each $\alpha<\kappa$,
$\domi{g_\alpha}{\le}{g}\in\cU_\alpha$.
\ee
\end{thm}
\begin{proof}
$1\Impl 2$: For each $\alpha<\kappa$ let
$Y_\alpha = \{f\in\NN : \domi{f}{<}{g_\alpha}\in\cU_\alpha\}$.
Then each $Y_\alpha\in\Dfin$, thus by (1) $Y=\Union_{\alpha<\kappa}Y_\alpha$
is not dominating. Let $g\in\NN$ be a witness for that.
In particular, for each $\alpha$ $g\nin Y_\alpha$, that is,
$\domi{g}{<}{g_\alpha}\nin\cU_\alpha$. As $\cU_\alpha$ is an ultrafilter,
we have that
$\domi{g_\alpha}{\le}{g}=\N\sm\domi{g}{<}{g_\alpha}\in\cU_\alpha$.

$2\Impl 1$: Assume that $Y=\Union_{\alpha<\kappa}Y_\alpha$ where
each $Y_\alpha\in\Dfin$. For each $\alpha$,
let $\cU_\alpha$ be an ultrafilter such that
$Y_\alpha/\cU_\alpha$ is bounded, say by $g_\alpha\in\NN$
\cite{CBC}.
By (2) let $g\in\NN$ be such that for each $\alpha<\kappa$,
$\domi{g_\alpha}{\le}{g}\in\cU_\alpha$.
Then $g$ witnesses that $Y$ is not dominating:
For each $f\in Y$, let $\alpha$ be such that $f\in Y_\alpha$.
Then $\domi{f}{\le}{g_\alpha}\in\cU_\alpha$, thus
$\domi{f}{<}{g}\spst\domi{f}{<}{g_\alpha}\cap
\domi{g_\alpha}{\le}{g}\in\cU_\alpha$; therefore
$\domi{f}{<}{g}$ is infinite.
\end{proof}

\begin{cor}\label{cofredprod}
Assume that $\cU$ is a nonprincipal ultrafilter on $\N$.
Then $\add(\Dfin,\fD)\le\cof(\NN/\cU)$.
\end{cor}
\begin{proof}
Assume that $\kappa<\add(\Dfin,\fD)$ and
let $\<g_\alpha : \alpha<\kappa\>$ be any $\kappa$-sequence
of elements of $\NN$.
For each $\alpha$ set $\cU_\alpha=\cU$.
Then by Theorem \ref{filteraddDfinD} there exists $g\in\NN$
such that for each $\alpha$, $\domi{g_\alpha}{\le}{g}\in\cU_\alpha=\cU$.
Thus $\<g_\alpha : \alpha<\kappa\>$ is not cofinal in $\NN/\cU$.
\end{proof}
%Corollary \ref{cofredprod} strengthens Mildenberger's result that
%for each nonprincipal ultrafilter $\cU$,
%$\g\le\cof(\NN/\cU)$ \cite{Mildenberger}.

\begin{cor}\label{cfd_bound}
$\add(\Dfin,\fD)\le\cf(\d)$.
\end{cor}
\begin{proof}
Canjar \cite{Canj} proved that there exists a nonprincipal ultrafilter
$\cU$ with $\cof(\NN/\cU)=\cf(\d)$. Now use Corollary \ref{cofredprod}.
\end{proof}

\begin{lem}\label{qo}
$g\in\NN$ avoids middles in $Y$ if, and only if,
for each $f\in Y$ $\domi{f}{<}{g}$ is infinite, and
the family $\{\domi{f}{<}{g} : f\in Y\}$
is linearly quasiordered by $\as$.
\end{lem}
%\begin{proof}
%This follows from the fact that for any $f,g,h\in\NN$,
%$\domi{f}{<}{g}\as\domi{h}{<}{g}$ if, and only if, $\mid{<}{\le}{g}{f}{h}$
%is finite.
%\end{proof}

\begin{thm}\label{addXD}
For any cardinal $\kappa$, the following are equivalent:
\be
\i $\kappa<\add(\fX,\fD)$;
\i For each $\kappa$-sequence $\<(g_\alpha,\cF_\alpha) :
\alpha<\kappa\>$,
such that each $g_\alpha\in\NN$, and for each $\alpha$
the restriction $\cF_\alpha\|\domi{0}{<}{g_\alpha}$
is large and linearly quasiordered by $\as$,
%if for each $\alpha$:
%\be
%\i $g_\alpha\in\NN$, % is not eventually equal to $0$,  --- follows from (c)
%\i $\cF_\alpha\sbst\inf$ is linearly quasiordered by $\as$; and
%\i The restriction $\cF_\alpha\|\domi{0}{<}{g_\alpha}$ is large.
%\ee
%Then
there exists $h\in\NN$ such that for each $\alpha<\kappa$,
the restriction $\cF_\alpha\|\domi{g_\alpha}{\le}{h}$ is large.
\ee
\end{thm}
\begin{proof}
$2\Impl 1$:
Assume that $Y=\Union_{\alpha<\kappa}Y_\alpha$ where each
$Y_\alpha\in\fX$.
For each $\alpha$ let $g_\alpha\in\NN$ be a function avoiding middles
in $Y_\alpha$, and
set $\cF_\alpha=\{\domi{f}{<}{g_\alpha} : f\in Y_\alpha\}$.
By Lemma \ref{qo}, $\cF_\alpha\sbst\inf$ is linearly quasiordered
by $\as$. As
$\cF_\alpha\|\domi{0}{<}{g_\alpha}=\cF_\alpha$,
the restriction is large and linearly quasiordered by $\as$.
By the assumption (2), there exists $h\in\NN$ such that for each $\alpha<\kappa$ and
each $f\in Y_\alpha$,
$\domi{f}{<}{g_\alpha}\cap\domi{g_\alpha}{\le}{h}$ is infinite; therefore
$h\not\le^* f$. Thus $h$ witnesses that $Y\in\fD$.

$1\Impl 2$: Replacing each $\cF_\alpha$ with $\cF_\alpha\|\domi{0}{<}{g_\alpha}$,
we may assume that each $A\in\cF_\alpha$
is an infinite subset of $\domi{0}{<}{g_\alpha}$.

For each $\alpha<\kappa$ let
$$Y_\alpha=\{f\in\NN : \domi{f}{<}{g_\alpha}\in\cF_\alpha %\|\domi{0}{<}{g_\alpha}
\}.$$
For each $A\in\cF_\alpha$ and each $h\in\NN$, define
\begin{equation}\label{genericf}
\tilde h(n) = \begin{cases}
g_\alpha(n)-1 & n\in A\\ %\cap\domi{0}{<}{g_\alpha}\\
\max\{g_\alpha(n),h(n)\} & \mbox{otherwise}
\end{cases}
\end{equation}
Then $\domi{\tilde h}{<}{g_\alpha}=A$,
%\cap\domi{0}{<}{g_\alpha}$,
%and the assumption (c) implies that this set
%is infinite.
%Moreover,
and $\domi{\tilde h}{<}{h}\sbst A$. %\cap\domi{0}{<}{g_\alpha}$.
Thus, for each $\alpha$,
$$\cF_\alpha = %\|\domi{0}{<}{g_\alpha}=
\{\domi{h}{<}{g_\alpha} : h\in Y_\alpha\}\sbst\inf.$$
As $\cF_\alpha$ is linearly quasiordered by $\as$, we
have by Lemma \ref{qo} that $g_\alpha$ avoids middles in $Y_\alpha$.
By (1), $Y=\Union_{\alpha<\kappa}Y_\alpha$ is not dominating; let
$h\in\NN$ be a witness for that.

For each $\alpha<\kappa$ and $A\in\cF_\alpha$,
let $\tilde h\in Y_\alpha$ be the function defined in
Equation \ref{genericf}.
Then $\tilde h\in Y$, therefore $\domi{\tilde h}{<}{h}$ is infinite.
By the definition of $\tilde h$,
$\domi{\tilde h}{<}{h}\sbst A\cap \domi{g_\alpha}{\le}{h}$;
therefore the restriction $\cF_\alpha\|\domi{g_\alpha}{\le}{h}$ is large.
\end{proof}

A nonprincipal ultrafilter $\cU$ is a \emph{simple $P_\kappa$ point}
if it is generated by a $\kappa$-sequence $\<A_\alpha : \alpha<\kappa\>\sbst\inf$
which is decreasing with respect to $\as$.
$\cU$ is a \emph{pseudo-$P_\kappa$ point} if every family $\cF\sbst\cU$
with $|\cF|<\kappa$ has a \pi{}.
Clearly every simple $P_\kappa$ point is a pseudo-$P_\kappa$ point.

\begin{cor}\label{cofP-point}
If $\cU$ is a simple $P_\kappa$ point, then $\add(\fX,\fD)\le\cof(\NN/\cU)$.
\end{cor}
\begin{proof}
Assume that $\lambda<\add(\fX,\fD)$.
Let $\<A_\beta : \beta<\kappa\>\sbst\inf$ be a $\kappa$-sequence which generates
$\cU$ and is linearly quasiordered by $\as$,
and set $\cF_\alpha = \cF = \{A_\beta : \beta<\kappa\}$ for all $\alpha<\lambda$.
Assume that $g_\alpha\in\NN$, $\alpha<\lambda$, are given. We will show that
these functions $g_\alpha$ are not cofinal in
$\NN/\cU$.

We may assume that for each $\alpha<\lambda$, $\domi{0}{<}{g_\alpha}=\N$.
Use Theorem \ref{addXD}
to obtain a function $h\in\NN$ such that for each $\alpha<\lambda$,
the restriction $\cF\| \domi{g_\alpha}{\le}{h}$ is large.
Assume that for some $\alpha<\lambda$, $\domi{g_\alpha}{\le}{h}\nin\cU$.
Then $\domi{h}{<}{g_\alpha}\in\cU$, thus there exists $\beta<\kappa$
such that $A_\beta\as \domi{h}{<}{g_\alpha}$, therefore
$A_\beta\cap \domi{g_\alpha}{\le}{h}$ is finite, a contradiction.
Thus $h+1$ witnesses that the functions $g_\alpha$ are not cofinal in
$\NN/\cU$, therefore $\lambda<\cof(\NN/\cU)$.
\end{proof}

In the remaining part of the paper we will consider the remaining standard cardinal
characteristics of the continuum (see \cite{HBK}).
Let $\u$ denote the minimal size of an ultrafilter base.
\begin{thm}\label{s_large}
It is consistent (relative to $ZFC$) that the following holds:
$$\u=\add(\Dfin,\fD)=\add(\fX,\fD)=\aleph_1<\aleph_2=\s=\mathfrak{c}.$$
Thus, it is not provable that $\s\le\add(\fX,\fD)$.
\end{thm}
\begin{proof}
In \cite{ShBl} a model of set theory is constructed where $\mathfrak{c}=\aleph_2$ and
there exist a simple $P_{\aleph_1}$ point and a simple $P_{\aleph_2}$ point.
The simple $P_{\aleph_1}$ point is generated by $\aleph_1$ many sets, thus
$\u=\aleph_1$. As $\b\le\u$, $\b=\aleph_1$ as well.

Nyikos proved that if there exists a pseudo $P_\kappa$ point $\cU$ and $\kappa>\b$, then
$\cof(\NN/\cU)=\b$ (see \cite{BlassMil}).
Thus by Corollary \ref{cofP-point}, $\add(\fX,\fD)\le\b=\aleph_1$ in this model.
In \cite{BlassMil} it is proved that if there exists a pseudo $P_\kappa$ point $\cU$,
then $\s\ge\kappa$. Therefore $\s\ge\aleph_2$ in this model.
\end{proof}

$\depth(\inf)$ is defined as the minimal cardinal $\kappa$ such that there exists
no $\subset^*$-decreasing $\kappa$-sequence in $\inf$.
(Thus, e.g., $\t<\depth(\inf)$.)
Each linearly quasiordered family $\cF\sbst\inf$ has a cofinal subfamily
which forms a $\subset^*$-decreasing sequence of length $<\depth(\inf)$.

\begin{thm}\label{ddepth}
~\be
\i If $\depth(\inf)<\d$, then $\add(\fX,\fD)=\d$.
\i If $\depth(\inf)=\d$, then $\cf(\d)\le\add(\fX,\fD)$.
\ee
\end{thm}
\begin{proof}
To prove (1) it is enough to show that for each $\kappa$ satisfying
$\depth(\inf)\le\kappa<\d$, we have that $\kappa<\add(\fX,\fD)$.
To prove (2) we will show that for each $\kappa<\cf(\d)$,
$\kappa<\add(\fX,\fD)$.
We will use Theorem \ref{addXD}, and prove both cases
simultaneously.

Assume that $\depth(\inf)\le\kappa<\d$ (respectively, $\kappa<\cf(\d)$).
Consider any
$\kappa$-sequence $\<(g_\alpha,\cF_\alpha) : \alpha<\kappa\>$
where each $g_\alpha\in\NN$, each $\cF_\alpha\sbst\inf$
is linearly quasiordered by $\as$, and the restriction
$\cF_\alpha\|\domi{0}{<}{g_\alpha}$ is large.
We must show that there exists $h\in\NN$ such that for each $\alpha<\kappa$,
the restriction $\cF_\alpha\|\domi{g_\alpha}{<}{h}$ is large.

Use the fact that $\depth(\inf)\le\kappa$ (respectively, $\depth(\inf)=\d$)
to choose for each $\alpha<\kappa$ a cofinal subfamily
$\tilde\cF_\alpha$ of $\cF_\alpha$ such that $|\tilde\cF_\alpha|<\kappa$
(respectively, $|\tilde\cF_\alpha|<\d$).

We may assume that each $g_\alpha$ is increasing.
For each $\alpha$ and each $A\in\cF_\alpha$,
let $\vec{A}\in\NN$ be the increasing enumeration of $A$.
The collection $\{g_\alpha\circ\vec{A} : \alpha<\kappa,\ A\in\cF_\alpha\}$
has less than $\d$ many elements and therefore cannot be dominating.
Let $h\in\NN$ be a witness for that.
Fix $\alpha<\kappa$. For all $A\in\cF_\alpha$, there exist
infinitely many $n$ such that
$$g_\alpha(\vec{A}(n))=g_\alpha\circ\vec{A}(n)<h(n)\le h(\vec{A}(n)),$$
that is, $A\cap\domi{g_\alpha}{<}{h}$ is infinite.
\end{proof}

%Let $\cov(\M)$ denote the minimal size of a family $\cF$ of
%meager (=first category) sets of reals such that $\cup\cF=\R$.

%\begin{cor}\label{covMdepth}
%If $\depth(\inf)<\cov(\M)$, then $\cov(\M)\le\add(\fX,\fD)$,
%and if $\depth(\inf)=\cov(\M)$, then $\cf(\cov(\M))\le\add(\fX,\fD)$.
%\end{cor}
%\begin{proof}
%If $\cov(\M)=\d$ then this is just Theorem \ref{ddepth}.
%Otherwise, $\cov(\M)<\d$ and we can use Theorem \ref{ddepth}(1)
%to get that in fact, $\add(\fX,\fD)=\d$.
%\end{proof}

\begin{thm}\label{cohenreals}
Assume that $V$ is a model of $CH$ and $\aleph_1<\kappa=\kappa^{\aleph_0}$.
%Assume that $\kappa$ is regular and
Let $\bbC_\kappa$ be the forcing notion which adjoins $\kappa$ many Cohen
reals to $V$.
Then in the Cohen model $V^{\bbC_\kappa}$, the following holds:
$$\add(\Dfin,\fD)=\s=\a=\non(\M)=\aleph_1<\cov(\M)=\add(\fX,\fD)=\mathfrak{c}.$$
\end{thm}
\begin{proof}
The assertions $\s=\a=\non(\M)=\aleph_1<\cov(\M)=\mathfrak{c}$
are well-known to hold in $V^{\bbC_\kappa}$, see \cite{HBK}.
It was proved by Kunen \cite{KunenPhD} that $V^{\bbC_\kappa}\models\depth(\inf)=\aleph_2$.
As $\cov(\M)\le\d$, we have that $\d=\mathfrak{c}=\kappa$ in this model.
If $\kappa=\aleph_2$, use Theorem \ref{ddepth}(1) and the fact that
$\d$ is regular in this model to obtain $\d\le\add(\fX,\fD)$.
Otherwise use Theorem \ref{ddepth}(2)
and the fact that $\depth(\inf)=\aleph_2<\kappa=\d$
to obtain this.

In \cite{Canj2, Roitman}
it is proved that there exists a nonprincipal ultrafilter $\cU$ in $V^{\bbC_\kappa}$
such that $\cof(\NN/\cU)=\aleph_1$. By Corollary \ref{cofredprod}, we have that
$\add(\Dfin,\fD)=\aleph_1$ in $V^{\bbC_\kappa}$.
\end{proof}

In particular, the cardinals $\add(\Dfin,\fD)$ and $\add(\fX,\fD)$ are not
provably equal.

\begin{cor}
It is not provable that $\add(\fX,\fD)\le\cf(\d)$.
\end{cor}
\begin{proof}
Use Theorem \ref{cohenreals} with $\kappa=\aleph_{\aleph_1}$.
In $V^{\bbC_\kappa}$, $\d=\mathfrak{c}=\aleph_{\aleph_1}$, therefore
$\cf(\d)=\aleph_1<\add(\fX,\fD)$ in this model.
\end{proof}

\begin{rem}
In the remaining canonical models of set theory which are used to distinguish
between the various cardinal characteristics of the continuum
(see \cite{HBK}),
$\max\{\b,\g\}=\d$ holds, and therefore $\add(\Dfin,\fD)=\add(\fX,\fD)=\d$ too.
These models show that none of the following is provable:
$\min\{\cov(\Null),\r\}\le\add(\fX,\fD)$ (\emph{Random reals} model),
$\add(\Dfin,\fD)\le\max\{\cov(\Null),\s\}$ (\emph{Hechler reals} model),
$\add(\Dfin,\fD)\le\max\{\non(\Null),\cov(\Null)\}$ (\emph{Laver reals} model),
and
$\add(\Dfin,\fD)\le\max\{\u,$ $\a,\non(\Null),\non(\M)\}$ (\emph{Miller reals} model).
\end{rem}

Collecting all of the consistency results,
we get that the only possible additional lower bounds on $\add(\fX,\fD)$ are
$\cov(\M)$ and $\e$ (observe that $\e\le\cov(\M)$ \cite{HBK}.)

\begin{prob}
Is $\cov(\M)\le\add(\fX,\fD)$?
And if not, is $\e\le\add(\fX,\fD)$?
\end{prob}

No additional cardinal characteristic can serve as an upper bound
on $\add(\Dfin,\fD)$.

Another question of interest is whether $\add(\Dfin,\fD)$ or $\add(\fX,\fD)$
appear in the lattice generated by the cardinal characteristics with
the operations of maximum and minimum. In particular, we have the following.

\begin{prob}\label{bgprob}
Is it provable that $\add(\Dfin,\fD)=\max\{\b,\g\}$?
\end{prob}

We have an indication that the answer to Problem \ref{bgprob}
is negative, but this is a delicate matter which will be treated
in a future work.

\iflogic
\else
\appendix
\section{Variations of the \emp{}}

\begin{defn}
For a subset $X$ of $\NN$, $g\in\NN$, and $R,S\in\{\le,<\}$, we
say that $g$ \emph{quasi avoids middles} in $X$ with respect to
$\<R,S\>$ if: \be \i $g$ is unbounded; \i for all $f,h\in X$ at
least one of the sets $\mid{R}{S}{g}{f}{h}$ and
$\mid{R}{S}{g}{h}{f}$ is finite. \ee A function $g\in\NN$
satisfying item (2) above is said to \emph{semi avoid middles} in
$X$ with respect to $\<R,S\>$. $X$ satisfies the \emph{quasi
excluded middle} property (respectively, \emph{semi excluded
middle} property) with respect to $\<R,S\>$ if there exists
$g\in\NN$ which quasi (respectively, semi) avoids middles in $X$
with respect to $\<R,S\>$. $\fx'_{R,S}$ (respectively,
$\fx''_{R,S}$) is the minimal size of a subset $X$ of $\NN$ which
does not satisfy the quasi (respectively, semi) excluded middle
property with respect to $\<R,S\>$.
\end{defn}

\begin{lem}\label{easytoo}
The following inequalities hold:
\be
\i $\fx'_{\le,\le}\le\fx'_{\le,<}\le\fx'_{<,\le}\le\fx'_{<,<}$,
\i $\fx''_{\le,\le}\le\fx''_{\le,<}\le\fx''_{<,\le}\le\fx''_{<,<}$;
\i For each $R,S\in\{\le,<\}$, $\fx_{R,S}\le\fx'_{R,S}\le\fx''_{R,S}$.
\ee
\end{lem}
\begin{proof}
(1) and (2) are proved as in Lemma \ref{easylemma}.
We will prove the first inequality of (3), the other one being immediate
from the definitions.
Assume that $Y\sbst\NN$ satisfies $|Y|<\fx_{R,S}$. Let $i\in\NN$ be
the identity function. Set $Y'=Y\cup\{i\}$. Then $|Y'|<\fx_{R,S}$,
thus there exists $g\in\NN$ which avoids middles in $Y'$.
In particular, the set $\domi{i}{R}{g}$ is infinite, thus $g$ is
unbounded, so $g$ quasi avoids middles in $Y$.
\end{proof}

\begin{thm}
The following assertions hold:
\be
\i Every subset $X$ of $\NN$ satisfies the \wemp{} with respect to
$\<<,\le\>$. Thus, $\fx''_{<,\le}=\fx''_{<,<}=\infty$.
\i $\fx'_{\le,\le}=\fx'_{\le,<}=\fx''_{\le,\le}=\fx''_{\le,<}=\b$;
\i $\fx'_{<,\le}=\fx'_{<,<}=\max\{\s,\b\}$.
\ee
\end{thm}
\begin{proof}
(1) Let $o\in\NN$ be the constant zero function.
Then for each $f,h\in\NN$, the set $\mid{<}{\le}{o}{f}{h}$
is finite.

(2) By Lemmas \ref{easylemma} and \ref{easytoo},
it is enough to show that $\fx''_{\le,<}\le\b$.
But this is, actually, what is shown in the proof of
Theorem \ref{equalsb}.

(3) By Lemmas \ref{easylemma}, \ref{s-le} and \ref{easytoo},
it is enough to show that $\fx'_{<,<}\le\max\{\s,\b\}$.
The proof is identical to the proof of Theorem \ref{sbthm}, since
for an unbounded $g\in\NN$, the set $[0<g]$ is
infinite, as required there.
\end{proof}
\fi

\begin{thebibliography}{00}

%+
\bibitem{JUBAR}
T.\ Bartoszy\'nski and H.\ Judah,
Set Theory: On the structure of the real line,
A.\ K.\ Peters, Massachusetts: 1995.

%+
\bibitem{HBK}
A.\ R.\ Blass,
\emph{Combinatorial cardinal characteristics of the continuum},
in: \textbf{Handbook of Set Theory} (M.\ Foreman, A.\ Kanamori, and M.\ Magidor, eds.),
Kluwer Academic Publishers, Dordrecht, to appear.

%+
\bibitem{BlassMil}
A.\ R.\ Blass and H.\ Mildenberger,
\emph{On the cofinality of ultrapowers},
Journal of Symbolic Logic \textbf{64} (1999),
727--736.

%+
\bibitem{ShBl}
A.\ R.\ Blass and S.\ Shelah,
\emph{There may be simple $P_{\aleph_1}$- and $P_{\aleph_2}$-points,
and the Rudin-Keisler ordering may be downward directed},
Annals of Pure and Applied Logic \textbf{33} (1987),
213--243.

%\bibitem{Canj1}
%R.\ M.\ Canjar,
%\emph{Model-Theoretic Properties of Countable Ultraproducts without the Continuum Hypothesis},
%Doctoral Dissertation, %Ph.D.\ Thesis,
%University of Michigan, 1982.

%+
\bibitem{Canj2}
R.\ M.\ Canjar,
\emph{Countable ultraproducts without CH},
Annals of Pure and Applied Logic \textbf{37} (1988), 1--79.

%+
\bibitem{Canj}
R.\ M.\ Canjar,
\emph{Cofinalities of countable ultraproducts: the existence theorem},
Notre Dame J.\ Formal Logic \textbf{30} (1989), 539--542.

%+
\bibitem{KaWe}
A.\ Kamburelis and B.\ W\c{e}glorz,
\emph{Splittings},
Archive for Mathematical Logic \textbf{35} (1996),
263--277.

%\bibitem{vD}
%E.\ K.\ van~Douwen,
%\emph{The  integers  and  topology},
%in: \textbf{Handbook of Set Theoretic Topology}
%(K.\ Kunen  and  J.\ Vaughan, Eds.),
%North-Holland, Amsterdam, 1984, 111--167.

%+
\bibitem{KunenPhD}
K.\ Kunen,
\emph{Inaccessibility Properties of Cardinals},
Doctoral Dissertation,
Stanford, 1968.

%+
\bibitem{Mildenberger}
H.\ Mildenberger,
\emph{Groupwise dense families},
Archive for Mathematical Logic \textbf{40} (2001),
93--112.

%+
\bibitem{Roitman}
J.\ Roitman,
\emph{Non-isomorphic H-fields from non-isomorphic ultrapowers},
Math.\ Z.\ \textbf{181} (1982),
93--96.

%+
\bibitem{ROTH2}
F.\ Rothberger,
\emph{On some problems of Hausdorff and of Sierpi\'nski},
Fund.\ Math.\ \textbf{35} (1948), 29--46.

%+
\bibitem{CBC}
M.\ Scheepers and B.\ Tsaban,
\emph{The combinatorics of Borel covers},
Topology and its Applications \textbf{121} (2002),
357--382.
\iflogic\else
\\ \arx{math.GN/0302322}
\fi

%\bibitem{Sh:207}
%S.\ Shelah,
%\emph{On cardinal invariants of the continuum},
%in: \textbf{Axiomatic set theory} (ed.\ Baumgartner, et.\ al.),
%Contemp.\ Mathematics \textbf{31} (1984),
%Amer. Math. Soc., Providence, RI,
%183--207.

%\bibitem{proper}
%S.\ Shelah,
%\emph{Proper and Improper Forcing} (2nd edition),
%Springer-Verlag, Berlin 1998.

%+
\bibitem{online}
S.\ Shelah and B.\ Tsaban,
\emph{Critical cardinalities and additivity properties of combinatorial notions of smallness}
(online version),
\arx{math.LO/0304019}

%+
\bibitem{tsaban1}
B.\ Tsaban,
\emph{A topological interpretation of $\mathfrak{t}$},
Real Analysis Exchange \textbf{25} (1999/2000), 391--404.
\iflogic\else
\\ \arx{math.LO/9705209}
\fi

%+
\bibitem{huremen1}
B.\ Tsaban,
\emph{A diagonalization property between Hurewicz and Menger},
Real Analysis Exchange \textbf{27} (2001/2002),
757--763.
\iflogic\else
\\ \arx{math.GN/0106085}
\fi

%+
\bibitem{tautau}
B.\ Tsaban,
\emph{Selection principles and the minimal tower problem},
Note di Matematica \textbf{22} (2003),
53--81.
\iflogic\else
\\ \arx{math.LO/0105045}
\fi

\newcommand{\Bc}[9]{\bibitem{#1} {#2}, \emph{#3}, in: \textbf{#4} (#5), #6 #7, #8--#9.}
\Bc{AddQuad}{B. Tsaban}{Additivity numbers of covering
properties}{Selection Principles and Covering Properties in
Topology} {L. Ko\v{c}inac, ed.}{Quaderni di Matematica 18, Seconda
Universita di Napoli, Caserta}{2006}{245}{282}

%\bibitem{Vaughan}
%J.\ Vaughan,
%\emph{Small uncountable cardinals and Topology},
%in: \textbf{Problems in Topology} (eds.\ Jan van Mill and G.M. Reed),
%North-Holland Pub.\ Co., Amsterdam: 1990,
%195--218.

\end{thebibliography}
\end{document}